\documentclass[a4paper]{article}

\usepackage{amsmath,amssymb,amsxtra}
\usepackage{mathtools}
\usepackage{mathrsfs}
\usepackage{latexsym}
\usepackage{bm}
\usepackage[amsmath,thref,thmmarks]{ntheorem}
\usepackage{etoolbox}
\usepackage{fancyhdr}
\usepackage{geometry}
\usepackage{titlesec}
\usepackage{tocloft}
\usepackage[
backend=biber,
style=alphabetic,
]{biblatex}
\usepackage{enumitem}
\usepackage{array}
\usepackage{tabularx}
\usepackage{tikz-cd}
\usepackage{graphicx}
\usepackage{authblk}
\definecolor{darkgreen}{rgb}{0,0.5,0}
\definecolor{darkblue}{rgb}{0,0,0.8}
\definecolor{darkred}{rgb}{0.8,0,0}
\definecolor{lightblue}{rgb}{0,0.6,0.8}
\usepackage[pdfencoding=auto,colorlinks,citecolor=darkgreen,linkcolor=darkblue,urlcolor=darkred]{hyperref}
\usepackage[capitalize]{cleveref}

\allowdisplaybreaks[2]
\newcommand\wan[1]{\widetilde{#1}}

\newcommand{\cl}[1]{\overline{#1}}

\newcommand{\bZ}{\mathbb{Z}}
\newcommand{\bQ}{\mathbb{Q}}
\newcommand{\bR}{\mathbb{R}}

\newcommand{\bF}{\mathbb{F}}

\newcommand\Spa[1]{\left\langle #1 \right\rangle}

\DeclareMathOperator{\im}{Im}

\DeclareMathOperator{\Aut}{Aut}

\DeclareMathOperator{\Gal}{Gal}

\DeclareMathOperator{\Sel}{Sel}
\DeclareMathOperator{\image}{im}

\DeclareMathOperator{\G}{G}
\DeclareMathOperator{\GL}{GL}
\titlespacing*{\section}{0pt}{20pt}{15pt}
\titlespacing*{\subsection}{0pt}{15pt}{12pt}
\titlespacing*{\subsubsection}{0pt}{12pt}{6pt}
{
\theorempreskip{5pt}
\theorempostskip{5pt}
\newtheorem{thm}{Theorem}[section]

\newtheorem{lemma}[thm]{Lemma}
\newtheorem{prop}[thm]{Proposition}
}
{
\theorempreskip{5pt}
\theorempostskip{5pt}
\theorembodyfont{\normalfont}
\newtheorem{defi}{Definition}[section]

\newtheorem{exam}{Example}[section]
\newtheorem{quest}{Question}[section]
}
{
\theorempreskip{5pt}
\theorempostskip{5pt}
\theorembodyfont{\normalfont}
\theoremsymbol{\ensuremath{\square}}
\newtheorem*{Proof}{Proof.}
}
{
\theoremstyle{break}
\theorempreskip{12pt}
\theorempostskip{5pt}
\theorembodyfont{\normalfont}
\theoremsymbol{\ensuremath{\square}}

}
\setlist[enumerate]{itemsep=0pt,parsep=2pt,topsep=2pt}
\setlist[itemize]{itemsep=3pt,parsep=3pt,topsep=5pt}
\makeatletter
\renewcommand{\@biblabel}[1]{[#1]\hfill}
\makeatother
  \DeclareFontFamily{U}{wncy}{}
    \DeclareFontShape{U}{wncy}{m}{n}{<->wncyr10}{}
    \DeclareSymbolFont{mcy}{U}{wncy}{m}{n}
    \DeclareMathSymbol{\Sh}{\mathord}{mcy}{"58} 
\fontsize{12pt}{20pt}\selectfont
\title{Computing the Iwasawa $\mu$-Invariants for Elliptic Curves over $\bZ_2$-Extensions}
\author[1]{Zichao Lin}
\author[2]{Mulun Yin}
\affil[1]{The University of Massachusetts Amherst}
\affil[2]{Morningside Center of Mathematics; Academy of Mathematics and Systems Sciences, Chinese Academy of Science}

\begin{document}
\maketitle
\abstract{Let $E$ be an elliptic curve defined over rational numbers. Further assume $E$ has good ordinary reduction at $p=2$ and $E[4]$ is reducible as a $G_\bQ$-representation. In this paper, we offer sufficient and necessary computational criteria for the algebraic Iwasawa $\mu$-invariant over the cyclotomic $\bZ_2$-extension of rational number to be 0, 1 or $\geq 2$. We also show that elliptic curves with isomorphic mod-4 representations have algebraic $\mu$-invariant equal to 0 if one of the curves does.}
\tableofcontents
\setlength{\abovedisplayskip}{5pt}
\setlength{\abovedisplayshortskip}{3pt}
\setlength{\belowdisplayskip}{5pt}
\setlength{\belowdisplayshortskip}{3pt}
\setlength{\jot}{3pt}
\section{Introduction}
\large
\subsection{Main results}
Let $E$ be an elliptic curve defined over $\bQ$. If $K$ is an algebraic extension of $\bQ$, fix a prime number $p$ where $E$ has good ordinary reduction. We can define the Selmer group $\Sel(E[p^\infty]/K)$, which fits into the following exact sequence \[0\to E(K)\otimes \bQ_p/\bZ_p\to \Sel(E[p^\infty]/K)\to \Sh(E/K)[p^\infty]\to 0,\] where $\Sh(E/K)$ is the Tate-Shafarevich group for $E$ over $K$. Using this exact sequence, we are able to study the rank of $E(K)$ and the order of $\Sh(E/K)$ using information about $\Sel(E[p^\infty]/K)$. \par
Now let $\bQ_\infty$ be the unique algebraic extension in a fixed algebraic closure with $\Gamma\coloneq\Gal(\bQ_\infty/\bQ)\simeq \bZ_p$. Set $\Lambda=\bZ_p[[\Gamma]]$, which is usually identified with the formal power series ring $\bZ_p[[T]]$ when a topological generator $\gamma$ of $\Gamma$ is sent to $1+T$, then we can view $\Sel(E[p^\infty]/\bQ_\infty)^\vee$ as a $\Lambda$-module. where $\Sel(E[p^\infty]/\bQ_\infty)^\vee$ denotes the Pontryagin dual of the $p^\infty$-Selmer group. According to the main result in \cite{Ka04}, $\Sel(E[p^\infty]/\bQ_\infty)$ is finitely generated $\Lambda$-cotorsion. Then by the fundamental structure theorem of finitely generated $\Lambda$-modules (see theorem 13.12 in \cite{Wa97}), there is a $\Lambda$-homomorphism from $\Sel(E[p^\infty]/\bQ_\infty)^\vee$ to \[(\bigoplus_{j=1}^t \Lambda/(f_i^{m_i}))\oplus (\bigoplus_{i=1}^s \Lambda/(p^{n_i}))\] with finite kernel and cokernel, with $f_i$ some irreducible polynomials in $\Lambda$. In this article, we will mainly focus on $\sum_{i=1}^s n_i$, which is called the algebraic $\mu$-invariant of $\Sel(E[p^\infty]/\bQ_\infty)$ and is denoted by $\mu_{E,p}$. \par
In \cite{GV00}, Greenberg and Vatsal proved the following two theorems related to the algebraic $\mu$-invariants of elliptic curves.
\begin{thm}\label{thm1: GV 1.4}
Assume $E$ and $E'$ are elliptic curves defined over $\bQ$, $p$ is an odd prime where both $E$ and $E'$ have good ordinary reduction. Assume $E[p]\simeq E'[p]$ as $G_\bQ$-modules, and they are irreducible as $\bF_p$-representations. Then $\mu_{E,p}=0$ if and only if $\mu_{E',p}=0$.
\end{thm}
This is theorem 1.4 of \cite{GV00}. \par
The first main result of this article is a generalization of theorem \ref{thm1: GV 1.4} to the case $p=2$. For simplicity, we will write $\mu_E$ for $\mu_{E,2}$. We have the following analogue. 
\begin{thm}\label{thm1: mu=0}
Let $E,E'$ be two elliptic curves defined over $\bQ$ with good ordinary reduction at $2$. If $E[4]\cong E'[4]$ as $G_{\bQ}$-representations, then $\mu_E=0$ if and only if $\mu_{E'}=0$.
\end{thm}
Compared to theorem \ref{thm1: GV 1.4}, we remove the irreducible condition and replace the naively expected $E[2]$ by $E[4]$. In fact, isomorphism between mod-$4$ representations is the best assumption we can have, as $E[2]\simeq E'[2]$ doesn't imply the equivalence between $\mu_{E,2}=0$ and $\mu_{E',2}=0$. See section \ref{sec 3} for an explicit example. \par
Next we consider the case where $E[p]$ is reducible as an $\bF_p$-representation with $\Phi$ a $G_\bQ$-invariant submodule. For a prime number $p$, let $I_p$ be the inertia group of any prime in $\cl{\bQ}$ above $p$. We say $\Phi$ is unramified at $p$ if $I_p$ acts on $A$ trivially and ramified otherwise. It is independent of the choice of the prime in $\cl{\bQ}$. We say $\Phi$ is odd if $\Phi$ lies inside the divisible subgroup of $E[p^\infty]$ such that complex conjugations act by $-1$ (in fact, it is independent of the choice of the complex conjugation). Otherwise, we say $\Phi$ is even. \par
In \cite{Gr99}, Greenberg proved the following two results on the algebraic $\mu$-invariants. The first provides a sufficient condition for $\mu_{E,p}=0$.
\begin{thm}\label{thm1: GV 1.3}
Assume $E$ is an elliptic curve defined over $\bQ$, and that $p$ is a prime where $E$ has good ordinary reduction. Assume $E$ has a $\bQ$-isogeny of degree $p$ with kernel $\Phi$. Suppose the action of $G_\bQ$ on $\Phi$ is either ramified at $p$ and even, or unramified at $p$ and odd. Then $\mu_{E,p}=0$.
\end{thm}
This is a combination of proposition 5.10 and 5.14 of \cite{Gr99}. See also theorem 1.3 in \cite{GV00}. In fact, the condition that $E$ has a $\bQ$-isogeny of degree $p$ is equivalent to the fact that $E[p]$ is reducible as a $G_\bQ$-representation over $\bF_p$. \par
In the situation where $\mu_{E,p}$ is positive, Greenberg proved the following result in \cite{Gr99}, which gives a lower bound of $\mu$-invariants.
\begin{thm}\label{thm1: Greenberg 5.7}
Assume $E$ is an elliptic curve defined over $\bQ$ with good ordinary reduction at $p$. Assume $E[p^\infty]$ contains a cyclic $G_\bQ$-invariant subgroup $\Phi$ of order $p^m$ which is ramified at $p$ and odd. Then $\mu_{E,p}\geq m$.
\end{thm}
This is a combination of propositions 5.7 and 5.13 of \cite{Gr99}. \par
In this article, we will focus on the cases where $p=2$ and $E[4]$ is reducible. Let $C$ be the intersection of the $p$-primary group of the formal group of elliptic curves over $\bQ_2$ and the divisible subgroup of $E[2^\infty]$ where complex conjugations act by $-1$. Under these assumptions and following the results in \cite{LY26}, we will offer a criterion for using the quadratic and biquadratic subextension of $\bQ(E[4])$, and possibly the order of $C$, to compute whether $\mu_{E,2}$ is $0$, $1$ or $\geq 2$, which will help us to pin down which elliptic curves in the isogeny class have $\mu$-invariant $0$ or $1$. Combining the results in \cite{LY26}, we are able to determine the algebraic $\mu$-invariants for all the elliptic curves in the isogeny class in an effective way.
\begin{thm}\label{thm1: E[4] red mu=0}
Assume $E$ is an elliptic curve defined over $\bQ$ with good ordinary reduction at $p=2$. Further assume $E[4]$ is reducible. There is a sufficient and necessary condition on whether $\mu_{E,2}=0$ if we know the quadratic subfields and biquadratic subfields in $\bQ(E[4])$.
\end{thm}
This result is seperated into several cases in section \ref{sec 4}. See theorems \ref{thm4: deg4 4*2 ram}, \ref{thm4: deg8 4}, \ref{thm4: deg4 2*2 Type B}, \ref{thm4: deg4 2*2 Type C}, \ref{thm4: deg8 2*2} and \ref{thm4: deg16 2} for details.
\begin{thm}\label{thm1: E[4] red mu=1}
Assume $E$ is an elliptic curve defined over $\bQ$ with good ordinary reduction at $p=2$. Further assume $E[4]$ is reducible. There is a sufficient and necessary condition on whether $\mu_{E,2}=1$ if we knows the quadratic subfields and biquadratic subfields in $\bQ(E[4])$ and the order of $C$.
\end{thm}
See theorems \ref{thm4: deg2 4*2}, \ref{thm4: deg4 4*2 unram}, \ref{thm4: deg4 2*2 Type A}, \ref{thm4: deg4 2*2 Type B}, \ref{thm4: deg4 2*2 Type C}, \ref{thm4: deg8 2*2}, \ref{thm4: deg8 2} and \ref{thm4: deg16 2} for details. \par
In fact, the observation of the relevance of the order of $C$ is from Greenberg. The condition in theorem \ref{thm1: GV 1.3} just means that $C$ has order $1$. 
\subsection{Strategy of the proof}
The proof has four main ingredients. The first one is proposition 2.4 and a generalization of proposition 2.8 of \cite{GV00}, which enable us to ignore all the primes other than 2 and the infinite prime, and compare the $\mu$-invariants of $\Sel(E[2^\infty]/\bQ_\infty)[4]$ with $\Sel(E[4]/\bQ_\infty)$. More precisely, we have the following result.
\begin{prop}\label{prop1: inside outside}
The Selmer groups $S(E[2^\infty]/\bQ_{\infty})^\vee/(2^n)$ and $S(E[2^n]/\bQ_\infty)^\vee$ have the same $\mu$-invariants.
\end{prop} 
This is proposition \ref{prop3: inside outside}.\par
From now on, we just need to compute the $\mu$-invariant of $\Sel(E[4]/\bQ_\infty)$ in order to study $\mu_{E,2}$.\par
The second ingradient is theorem \ref{thm1: GV 1.3} and \ref{thm1: Greenberg 5.7}, which enable us to give sufficient conditions for $\mu_{E,2}=0$ and $\mu_{E,2}\geq 1$. The philosophy of these two results is that the $G_\bQ$-invariant subgroup $E[4]$ can contribute to $\mu_{E,2}$. In other words, understanding how $\Gal(\bQ(E[4])/\bQ)$ acts on $E[4]$, especially the action on the $G_\bQ$-invariant part is helpful for computing the $\mu$-invariant. \par
The third ingredient is Serre's results in \cite{Se72}. More precisely, $I_2$ acts on the ramified part via the cyclotomic character, and acts on the quotient of the ramified part by the trivial character. Moreover, the ramified part always forms a sub-representation of $G_{\bQ_2}$.\par
For the infinite prime, since the cyclotomic character is odd, there is always a subrepresentation of $E[4]$ where complex conjugations act by $-1$, and a quotient representation where complex conjugations act by $1$. \par
The last ingradient is the Greenberg's conjecture proved in \cite{LY26}. In fact, we need an upper bound of the $\mu$-invariants in order to compute the precise $\mu$-invariants when $E[4]$ is reducible. This is accomplished by applying proposition 3.2 and theorem 5.5 in \cite{LY26}.
\subsection{Relations to previous works}
The study of the $\mu$-invariant of $\Sel(E[p^\infty]/\bQ_\infty)$ starts from Greenberg \cite{Gr99}. For odd $p$, in \cite{GV00}, Greenberg and Vatsal also computed similar results for the analytic $\mu$-invariant, which is defined using Mazur and Swinnerton-Dyer's construction of the p-adic L-function of the Hasse-Weil L-function in \cite{MSD74}. Using Greenberg's criterion of the lower bound of the $\mu$-invariants, Kramer constructed some explicit families of elliptic curves defined over $\bQ$ with $\mu_{E,2}\geq n$, where $n=1,2,3,4$ in \cite{Kr99}. \par
Under the assumption that $E$ has good ordinary reduction at $p$ and $E[p]$ is reducible, a typical case left untouched by Greenberg and Vatsal is when $E[p]$ can be placed in the nonsplit exact sequence \[0\to \bZ/p\bZ\to E[p]\to \mu_p\to 0.\] Under this assumption, Trifkovic found infinitely elliptic curves with $\mu_{E,p}=0$ where $p=3$ or 5 in \cite{Tr05} using the global duality pairing for finite flat group scheme. Using a different method, in \cite{Ha04} Hachimori proved that $\mu_{E,3}=0$ is equivalent to the vanishing of the $\mu$-invariant of a Galois group with restricted ramification over $K_\infty$, where $K$ is the real cubic subfield of $\bQ(E[3])$. Recently, in \cite{LY26}, the authors are able to prove the Greenberg's conjecture for $p=2$ under the assumption that $E[2]$ is reducible in the good ordinary case. \par
This article contains two novelties. Compared to the previous works which are all considering mod-$p$ representation, this is the first one that considers mod-$p^2$ representation. The reason lies in the special property of $p=2$. Let $C_p$ denote the kernel of the reduction map \[E[p^\infty]\to \wan{E}[p^\infty],\] then $G_{\bQ_p}$ acts on $C_p$ via the cyclotomic chracter $\omega_p$. But $\omega_2$ is trivial on $C_2[2]$, which doesn't provide much information. This is the reason why we need to use the mod-4 representation to get more information about $C_2$.  \par
The second novelty is that we are able obtain an upper bound of the $\mu$-invariants even in the case where $\mu_{E,2}$ is positive, which enable us to provide explicit examples with $\mu_{E,2}=1$. The new input is corollary of the main result in \cite{LY26}. In \cite{LY26}, we are able to show that if we assume $E[2]$ is reducible, then all the $\mu$-invariants are from the odd and ramified submodule. So in order to give an upper bound of the $\mu$-invariant, we just need to offer criterions to bound the order of the intersection of the odd submodule and ramifies submodule from above.
\subsection{Further developments}
The first question one can ask is whether the results in this article can be generalized to the case $p\geq 3$. In \cite{LY26b}, the authors are able to prove Greenberg's conjecture at all primes $p\geq 3$ with ordinary reduction under the assumption that $E[p]$ is reducible. In order to obtain a computational result as in this article, one needs a similar result for $p\geq3$ on how to characterize the odd and ramified submodules using subextensions of $\bQ(E[p])$. Since still $E[p]$ may well determine the $\mu$-invariant only up to 1, probably one also needs to analyze the mod-$p^2$ representation in order to obtain a concrete description of when $\mu=1$. \par
The second question is whether one is able to compute the exact value of the $\mu$-invariant at $p=2$ when $\mu_{E,2}\geq 2$. In fact, our work can be viewed as evidence that $E[4]$ can control the $\mu$-invariant up to $1$, so in order to find examples with $\mu_{E,2}=2$, one would naturally need to analyze $E[8]$. In other words, one needs to compute the mod-$8$ representation. Probably this can be done with some assistance from code programming using Magma or PARI. \par
The third question is whether one can use only the information of the Galois representations to compute the orders of the intersection of the odd and ramified submodules. According to the concrete examples in section \ref{5.3}, it seems that one also needs at least $E[8]$ to give a complete solution. \par
If $E$ has multiplicative reduction at $p=2$, it is also known that $\Sel(E[p^\infty]/\bQ_\infty)^\vee$ is a finitely generated $\Lambda$-torsion module, so one is still able to define the $\mu$-invariant for it. It is natural to ask whether one can obtain a similar criterion if $E$ has multiplicative reduction at $p=2$.
\subsection{Organization of the paper}
In section \ref{sec 2}, we will recall the definitions of Selmer groups and Iwasawa invariants. The key result of that section is the relation between $2^\infty$-Selmer groups and $2^n$-Selmer groups. In section \ref{sec 3} we will prove a theorem that generalizes the results in \cite{GV00} on the algebraic $\mu$-invariants of elliptic curves with isomorphic mod-$p$ representations. Section \ref{sec 4} is the main section of this article. It is about computing the $\mu$-invariant of an elliptic curve $E$ with $E[4]$ reducible and good ordinary reduction at $p=2$. We will separate it into four cases based on the structure of $E(\bQ)[4]$. Section \ref{sec 5} is about numerical examples.

\subsection{Acknowledgements}
The authors would like to thank Robert Pollack, Tom Weston and David Zureick-Brown for helpful discussions and advices. This project was initiated when the authors were attending the Arithmetic Geometry workshop on Iwasawa theory at the University of North Texas in May 2025. The authors are grateful for UNT's organizing and hosting, and thank Jeff Hatley for proposing this question. This work is part of the first author's forthcoming Ph.D thesis.

\section{Selmer groups and Iwasawa invariants}\label{sec 2}
\large
Let $E/\bQ$ an elliptic curve with good ordinary reduction at $p=2$. In this section we will define relevant Selmer groups and Iwasawa invariants and prove some basic properties.

\subsection{Selmer groups}
Let $\bQ_\infty$ be the cyclotomic $\bZ_p$-extension of $\bQ$, and let $\eta$ be a place of $\bQ_\infty$. Consider the Kummer map \[\kappa: E(\bQ_\infty)\otimes \bQ_2/\bZ_2\to H^1(\bQ_\infty, E[2^\infty]), \] and its local at $\eta$ version \[\kappa_\eta : E((\bQ_\infty)_\eta)\otimes \bQ_2/\bZ_2\to H^1((\bQ_\infty)_\eta, E[2^\infty]).\] 
\begin{defi}
Define the 2-primary Selmer group to be \[\Sel(E[2^\infty]/\bQ_\infty)=\ker\Bigl(H^1(\bQ_\infty, E[2^\infty])\to \prod_\eta H^1((\bQ_\infty)_\eta, E[2^\infty])/\image(\kappa_\eta) \Bigr). \]
\end{defi} 
According to \cite{Gr89} and \cite{Gr99}, we can give a equivalent definition of the $2$-primary Selmer group. If $l\neq2$, we let \[H_l(\bQ_\infty, E[2^\infty])=\prod_{\eta\mid l} H^1((\bQ_\infty)_\eta, E[2^\infty]).\] Let $\eta_2$ be the unique prime of $\bQ_\infty$ above $\bQ$, fix a prime $\pi$ above $\eta_2$ in $\cl{\bQ}$ and let $I_2$ be the inertia subgroup with respect to $\pi$ in $G_{\bQ_\infty}$. Let $\wan{E}$ be the reduction of $E$ modulo $2$. Define \[L_2=\ker\Bigl(H^1((\bQ_\infty)_{\eta_2}, E[2^\infty])\to H^1(I_2, \wan{E}[2^\infty])\Bigr)\] and \[H_2=H^1((\bQ_\infty)_{\eta_2}, E[2^\infty])/L_2.\] Fix $\Sigma$ a finite set of primes containing $2$, $\infty$ and all the primes where $E$ has bad reduction. Denote by $\bQ_\Sigma$ the maximal extension of $\bQ$ unramified outside $\Sigma$.
\begin{prop}\label{prop2: Equivalent definition}
An equivalent definition of the $2$-primary Selmer group over $\bQ_\infty$ is \[\Sel(E[2^\infty]/\bQ_\infty)=\ker\Bigl(H^1(\bQ_\Sigma/\bQ_\infty, E[2^\infty])\to \prod_{l\in \Sigma} H_l(\bQ_\infty, E[2^\infty])\Bigr). \]
\end{prop}

This is in section 2 of \cite{Gr89}.

For further computational purpose, we also need to define the non-primitive $2$-primary Selmer group over $\bQ_\infty$ and the $2^n$-Selmer group. Let $\Sigma_0=\Sigma-\{ 2,\infty\}$. 
\begin{defi}
Define \[S(E[2^\infty]/\bQ_{\infty})=\Sel^{\Sigma_0}(E[2^\infty]/\bQ_\infty)=\ker\Bigl( H^1(\bQ_\Sigma/\bQ_\infty, E[2^\infty])\to \prod_{l\in \Sigma-\Sigma_0} H_l(\bQ_\infty)\Bigr).\]
\end{defi}

Define \[\kappa_n: E(\bQ_\infty)/2^nE(\bQ_\infty)\to H^1(\bQ_\infty, E[2^n])\]to be the $2^n$-Kummer map, with\[\kappa_{n,\eta}: E((\bQ_\infty)_\eta)/2^nE((\bQ_\infty)_\eta)\to H^1((\bQ_\infty)_\eta, E[2^n])\] the local at $\eta$ version.
\begin{defi}
Define the non-primitive $2^n$-Selmer group to be \[S(E[2^n]/\bQ_\infty)=\ker\Bigl(H^1(\bQ_\infty, E[2^n])\to \prod_{\eta\in \Sigma-\Sigma_0} H^1((\bQ_\infty)_\eta, E[2^n])/\image(\kappa_{n,\eta})\Bigr). \]
\end{defi}

\subsection{Iwasawa invariants}

Let $\Lambda=\bZ_p[[\Gal(\bQ_\infty/\bQ)]]=\bZ_p[[\Gamma]]$, which is usually identified with the formal power series ring $\bZ_p[[T]]$ when a topological generator $\gamma$ of $\Gamma$ is sent to $1+T$. Then the Selmer groups $\Sel(E[2^\infty]/\bQ_\infty), S(E[2^\infty]/\bQ_\infty), \Sel(E[2^n]/\bQ_\infty)$ have natural $\Lambda$-module structures. According to \cite{Ka04}, $\Sel(E[2^\infty]/\bQ_\infty)$ is a co-finitely generated $\Lambda$-cotorsion module. In other words, there exists a pseudo-isomorphism (that is, a $\Lambda$-morphism with finite kernel and cokernel) \[\Sel(E[2^\infty]/\bQ_\infty)^\vee\to \Bigl(\bigoplus_{i=1}^t \Lambda/(f_i)\Bigr)\bigoplus \Bigl(\bigoplus_{j=1}^s \Lambda/(2^{n_j})\Bigr),\]
where $(-)^\vee$ denotes Pontryagin dual. We define the $\mu$-invariant of the elliptic curve $E$ to be \[\mu_E=\sum_{j=1}^s n_j.\]

The first property of $\mu$-invariant is that it is the preserved under pseudo-isomorphism of torsion $\Lambda$-modules. 
\begin{lemma}\label{lem2: same mu}
Let $A,B$ be two $\Lambda$-torsion modules, suppose $f: A\to B$ is a pseudo-isomorphism, then $A$ and $B$ have the same $\mu$-invariants.
\end{lemma}
\begin{Proof}
Consider the following exact sequence of $\Lambda$-modules
\[\begin{tikzcd}
0 \ar{r} & F \ar{r} & A \ar{r}{f} & B \ar{r} & G \ar{r} & 0 \\
\end{tikzcd}\] with $F$ and $G$ finite. Let $\Lambda_2$ be the localization of $\Lambda$ at the prime $(2)$. Tensor the exact sequence by $\Lambda_2$. Since $A$ is finite, there exists an $n$ such that $2^nA=0$, hence $A\otimes \Lambda_2=0$. Similarly, $B\otimes \Lambda_2=0$. Therefore we are left with $A\otimes \Lambda_2\simeq B\otimes \Lambda_2$. The result follows from lemma 15.18 of \cite{Wa97}.
\end{Proof}

The reason why we introduce the non-primitive Selmer group is that it has the same $\mu$-invariant as the usual Selmer group, and is often more convenient to deal with.

\begin{prop}\label{prop2: same mu}
$S(E[2^\infty]/\bQ_{\infty})$ is $\Lambda$-cotorsion. The $\mu$-invariants of $S(E[2^\infty]/\bQ_{\infty})^\vee$ and $\Sel_2(E/\bQ_{\infty})^\vee$ are equal.
\end{prop}
\begin{Proof}
This is corollary 2.3 in \cite{GV00}. \par
\end{Proof}

\subsection{Some results from \cite{LY26}}
In this section, we will state some results from \cite{LY26} that are related to the $\mu$-invariants of Galois cohomology groups relevant to the Selmer groups. First we will give an alternative description of the local conditions defining $S(E[2^n]/\bQ_\infty)$. Denote \[C_2=\ker\Bigl(E[2^\infty]\to \wan{E}[2^\infty]\Bigr).\] Since $E$ has good ordinary reduction at $2$, $C_2$ has $\bZ_2$-corank $1$. Let $D_2=E[2^\infty]/C_2$, which is a quotient module of $E[2^\infty]$ with trivial $I_2$-action. For an infinite prime $v$ of $\bQ_\infty$, denote by $C_{\infty, v}$ the divisible part of the minus part of $E[2^\infty]$ with respect to the complex conjugation induced by $v$. Let $D_{\infty,v}=E[2^\infty]/C_{\infty,v}$. Finally let $\Delta$ be the discriminant of $E$.

\begin{prop}\label{prop2: description of E[4]}
If $\Delta>0$, then
\[S(E[2^n]/\bQ_\infty)\cong \ker\Bigl(H^1(\bQ_\Sigma/\bQ_\infty, E[2^n])\to H^1(I_2, D_2[2^n])\times \prod_{v\mid \infty} H^1(\bR, D_{\infty,v}[2^n])\Bigr).\]
If $\Delta<0$, then \[S(E[2^n]/\bQ_\infty)\cong \ker\Bigl(H^1(\bQ_\Sigma/\bQ_\infty, E[2^n])\to H^1(I_2, D_2[2^n])\Bigr).\]
\end{prop}

\begin{Proof}
This is proposition 2.5 in \cite{LY26}.
\end{Proof}

Next we will state a result about the $2^n$-division field of $E$, which will be used frequently in the classification of the $\mu$-invariants in section \ref{sec 4}.

\begin{lemma}\label{lem2: not Q(sqrt2)}
Assume $E$ an elliptic curved defined over $\bQ$ with good ordinary reduction at $p=2$. Suppose $K:=\bQ(E[2])$ is a quadratic field, then $K=\bQ(\sqrt{m})$ with $2\nmid m$.
\end{lemma}
\begin{Proof}
This follows from the proof of lemma 4.7 of \cite{LY26}. In fact the proof still works under the assumption that $K$ is imaginary quadratic.
\end{Proof}

\section{A Greenberg-Vatsal type result for $p=2$}\label{sec 3}
\large
In this section, we will give a generalization of the following result by Greenberg and Vatsal in \cite{GV00}.
\begin{thm}\label{thm3: GV}
Assume $E_1$ and $E_2$ are elliptic curves defined over $\bQ$. Let $p$ be an odd prime where both $E_1$ and $E_2$ have good ordinary reduction. Assume $E_1[p]\simeq E_2[p]$ as Galois modules and these are irreducible. Then $\mu_{E_1}=0$ if and only if $\mu_{E_2}=0$.
\end{thm}

This result is not true for $p=2$. For counter-examples, let $E_1=15A5$, $E_2=15A6$ in the LMFDB table. Since $E_1(\bQ)[2]\simeq E_2(\bQ)[2]\simeq (\bZ/2\bZ)^2$, their mod-2 representations are isomorphic. However, using the compuation from section \ref{sec 5}, we will see that $\mu_{E_1}=1$, but $\mu_{E_2}=0$. \\
In order to derive a similar result for $p=2$, we need to compare $E[4]$ instead. For simplicity, we denote $\Sel(E[2^\infty]/\bQ_\infty)$, $S(E[2^\infty]/\bQ_\infty)$ and $S(E[2^n]/\bQ_\infty)$ by $\Sel(E[2^\infty])$, $S(E[2^\infty])$ and $S(E[2^n])$ respectively. The main theorem of this section may be stated as follows.

\begin{thm}\label{thm3: mu=0}
Let $E,E'$ be two elliptic curves defined over $\bQ$ with good ordinary reduction at 2. If $E[4]\cong E'[4]$ as $G_{\bQ}$-representations, then $\mu_E=0$ if and only if $\mu_{E'}=0$.
\end{thm}

The key ingredient of the proof is proposition 2.6 in \cite{LY26}, which is a modification of proposition 2.8 in \cite{GV00}. Using this result, we are able to compare the $\mu$-invariants of $\bigl(\Sel([2^{\infty}])[2^n]\bigr)^\vee$ and $\bigl(S([2^n]))^\vee$. The main difference between our result and the one in \cite{GV00} is that we need to tackle the infinite primes for $p=2$, which do not show up in the Selmer groups for $p\geq 3$. 

\begin{prop}\label{prop3: inside outside}
The Selmer groups $S(E[2^\infty]/\bQ_{\infty})^\vee/(2^n)$ and $S(E[2^n]/\bQ_\infty)^\vee$ have the same $\mu$-invariants.
\end{prop} 

\begin{Proof}
See proposition 2.6 of \cite{LY26}.
\end{Proof}

Now we are able to prove theorem \ref{thm3: mu=0}.
\begin{Proof}[theorem \ref{thm3: mu=0}]
Combining proposition \ref{prop3: inside outside} with the structure theorem of finitely generated $\Lambda$-modules, one sees that $S(E[2]/\bQ_\infty)^\vee$ has $\mu$-invariant 0 if and only if $\mu_E=0$. Therefore, it suffices to show that $E[4]\cong E'[4]$ implies that the Selmer groups $S(E[2]/\bQ_\infty)$ and $S(E'[2]/\bQ_\infty)$ are isomorphic.\par
Let $f:E[4]\to E'[4]$ be a $G_\bQ$-isomorphism. Let $\{P,Q\}$ be a basis of $E[4]$ and set $P'=f(P)$, $Q'=f(Q)$. Since $E[4]\cong E'[4]$, we have $H^1(\bQ_\Sigma/\bQ_\infty, E[2])\cong H^1(\bQ_\Sigma/\bQ_\infty, E'[2])$. \par
To show that the Selmer groups coincide, it suffices to show that the the local conditions at $2$ and infinite primes agree. Since $E$ has good ordinary reduction at $2$, there exists some $\sigma\in I_2$ such that $\sigma$ acts on $E[4]$ by $\begin{pmatrix}
-1 & t \\
0 & 1
\end{pmatrix}$ with $0\leq t\leq 3$. If $t=1$ or $3$, then the only $4$-torsion on which $I_2$ acts via the cyclotomic character is $\Spa{P}$, hence it is fixed under $f$. If $t=0$ or $2$, $C_2[4]$ can be either $\Spa{P}$ or $\Spa{P+2Q}$, but in either case, $C_2[2]$ is $\Spa{2P}$, which is still fixed by $f$. For infinite primes, since complex conjugation acts on $C_\infty$ by $-1$, and it acts on $E[4]$ also by $\begin{pmatrix}
-1 & t \\
0 & 1
\end{pmatrix}$, we can use the same argument as for $I_2$.
\end{Proof}

\section{Computing $\mu_E$ for elliptic curves with reducible $E[4]$}\label{sec 4}
From now on, we assume $E/\bQ$ is an elliptic curve with good ordinary reduction at $p=2$. Furthermore, we assume $E[4]$ is reducible as a $\G_\bQ$-module. Our main result of this section is a full classification of the elliptic curves with $\mu_E=0$, $\mu_E=1$ and $\mu_E\geq 2$.
\subsection{Some preliminary results}
Let $K=\bQ(E[4])$ be the field of definition of the $4$-torsion points of $E$. Since $E[4]$ is reducible, we have $E(\bQ)[4]\neq 0$, and hence $E(\bQ)[4]\simeq \bZ/4\bZ\times \bZ/2\bZ,\bZ/4\bZ, \bZ/2\bZ\times \bZ/2\bZ$ or $\bZ/2\bZ$. We will consider these cases separately. The key techniques that we will apply are the following two results from \cite{LY26}.
\begin{prop}\label{prop4: delta}
Let $\alpha$ be an order 2 $G_\bQ$-submodule of $E[2^\infty]$ that fits into the following exact sequence \[0\to \alpha\to E\to E'.\] Then
\begin{enumerate}
\item If $\alpha$ is odd and ramified, then $\mu_E-\mu_{E'}=1$.
\item If $\alpha$ is not ramified and not odd, then $\mu_E-\mu_{E'}=-1$.
\item Otherwise, $\mu_E=\mu_{E'}$.
\end{enumerate}
\end{prop}
\begin{Proof}
See proposition 3.2 of \cite{LY26}.
\end{Proof}

\begin{prop}\label{prop4: mu=m}
Assume $E$ is defined over $\bQ$ and has good reduction at $p=2$. Further assume $E[2]$ is reducible. Then $\mu_E=n$ if and only if the largest $G_\bQ$-invariant ramified and odd subgroup $\Phi$ has order $2^n$.
\end{prop}
\begin{Proof}
This is theorem 5.5 of \cite{LY26}.
\end{Proof}

\subsection{Computing $\mu_E$ when $E(\bQ)[4]\simeq \bZ/4\bZ\times \bZ/2\bZ$}
Let $\{P,Q\}$ be a basis of $E[4]$. Let \[\rho:\Gal(\bQ(E[4])/\bQ)\to \Aut(E[4])\simeq \GL_2(\bZ/4\bZ).\] If $E(\bQ)[4]\simeq \bZ/4\bZ\times \bZ/2\bZ$, we may assume $\im(\rho)$ is of the form 
$\begin{pmatrix}
1 & b \\
0 & d
\end{pmatrix}$
with $b\in \{0,2\}, d\in \{1,3\}$. According to the Weil pairing, the determinant of the action should be the cyclotomic character, hence $d$ is not always $1$. Therefore, $[\bQ(E[4]):\bQ]=2$ or $4$.\par
When $[\bQ(E[4]):\bQ]=2$, using the Weil paring, we have $\bQ(E[4])=\bQ(i)$ and the nontrivial element of $\Gal(\bQ(i)/\bQ)$ has image $\begin{pmatrix}
1 & 0 \\
0 & 3
\end{pmatrix}$ under $\rho$. Since the inertia group acts on $C_2$ via the $4$-th cyclotomic character, $C_2[4]=\Spa{Q}$ or $\Spa{Q+2P}$, and by changing a basis if necessary, we may assume $C_2[4]=\Spa{Q}$. Also, since a complex conjugation $c_v$ acts on $C_{\infty,v}$ by $-1$, $C_{\infty,v}[4]=\Spa{Q}$ or $\Spa{Q+2P}$. Notice that both $\Spa{Q}$ and $\Spa{Q+2P}$ are $G_\bQ$-invariant, hence $C_{\infty,v}$ is independent of the $v$ in $\bQ_\infty$ and we may denote it simply as $C_\infty[4]$.

\begin{thm}\label{thm4: deg2 4*2}
Suppose $E(\bQ)[4]\simeq \bZ/4\bZ\times \bZ/2\bZ$ and $[\bQ(E[4]):\bQ]=2$. If the intersection of $C_2[4]$ and $C_\infty[4]$ has order 4, then $\mu_E\geq 2$. If the intersection has order 2, then $\mu_E=1$.
\end{thm}

\begin{Proof}
Both cases follow directly from proposition \ref{prop4: mu=m}. In fact, all the cyclic submodules of $E[4]$ are $G_\bQ$-invariant. If $C_2[4]\cap C_\infty[4]$ has order $4$, it means the maximal $G_\bQ$-invariant odd and ramified submodule has order $\geq 4$. Therefore, $\mu_E\geq 2$. If $C_2[4]\cap C_\infty[4]$ has order $2$, then $C_2[4]\cap C_\infty[4]$ is exactly the largest odd and ramified submodule of $E[2^\infty]$. Hence by proposition \ref{prop4: mu=m}, $\mu_E=1$.
\end{Proof}

Next we assume $[\bQ(E[4]):\bQ]=4$, then $\im(\rho)\simeq (\bZ/2\bZ)^2$ and is generated by $\begin{pmatrix}
1 & 0 \\
0 & 3
\end{pmatrix}$ and $\begin{pmatrix}
1 & 2 \\
0 & 1
\end{pmatrix}$. Here we still assume $C_2[4]=\Spa{Q}$. Hence $\bQ(E[4])=\bQ(\sqrt{m},i)$ with $m$ a square-free positive integer. We first prove that $2\nmid m$.
\begin{lemma}\label{lem4: e=2}
Under the assumptions listed above, $\bQ(E[4])=\bQ(\sqrt{m},i)$ with $2\nmid m$.
\end{lemma}

\begin{Proof}
It suffices to show that $\bQ(E[4])/\bQ$ is not totally ramified. Suppose not, then both $\begin{pmatrix}
1 & 0 \\
0 & 3
\end{pmatrix}$ and $\begin{pmatrix}
1 & 2 \\
0 & 3
\end{pmatrix}$ are in the images of $I_2$ under $\rho$. Since $I_2$ acts on the odd part via the cyclotomic character, assuming $\tau\in \Gal(\bQ(E[4])/\bQ)$ has image $\begin{pmatrix}
1 & 2 \\
0 & 3
\end{pmatrix}$ under $\rho$, then $\rho(\tau)(Q)=2P+3Q$ lies inside $\Spa{Q}$, contradiction arises.
\end{Proof}

Now we can state the main result under the assumptions $E(\bQ)[4]\simeq \bZ/4\bZ\times \bZ/2\bZ$ and $[\bQ(E[4]):\bQ]=4$. We will seperate it into two theorems depending on whether $2$ ramifies in $\bQ(\sqrt{m})$.

\begin{thm}\label{thm4: deg4 4*2 ram}
Assume $E(\bQ)[4]\simeq \bZ/4\bZ\times \bZ/2\bZ$ and $[\bQ(E[4]):\bQ]=4$. Let $K=\bQ(\sqrt{m})$ be the unique real quadratic subfield of $\bQ(E[4])$. If 2 is ramified in $K$, then $\mu_E=0$.
\end{thm}

\begin{Proof}
By proposition \ref{prop4: mu=m}, it suffices to show $2Q\notin C_{\infty}[2]$. By the Weil pairing, $\rho(\Gal(\bQ(E[4])/\bQ(i)))$ is generated by $\begin{pmatrix}
1 & 2 \\
0 & 1
\end{pmatrix}$. Since the $I_2$ action on $E[4]$ factors through $\bQ(E[4])/\bQ(\sqrt{-m})$, we can see that $\rho(\Gal(\bQ(E[4])/\bQ(\sqrt{-m})$ is generated by $\begin{pmatrix}
1 & 0 \\
0 & 3
\end{pmatrix}$. Therefore, the action of a complex conjugation on $E[4]$ factors through the nontrivial element in $\Gal(\bQ(E[4])/\bQ(\sqrt{m})$, which has image $\begin{pmatrix}
1 & 2 \\
0 & 3
\end{pmatrix}$ under $\rho$, hence $C_\infty[4]$ is either $\Spa{P+Q}$ or $\Spa{P+3Q}$. In either case, $C_\infty[2]=\Spa{2P+2Q}$, therefore $2Q\notin C_\infty[2]$.
\end{Proof}

\begin{thm}\label{thm4: deg4 4*2 unram}
Assume $E(\bQ)[4]\simeq \bZ/4\bZ\times \bZ/2\bZ$ and $[\bQ(E[4]):\bQ]=4$. Let $K=\bQ(\sqrt{m})$ be the unique real quadratic subfield of $\bQ(E[4])$. If 2 is unramified in $K$, then $\mu_E=1$.
\end{thm}

\begin{Proof}[theorem \ref{thm4: deg4 4*2 unram}]
We will use an argument close to the proof of theorem \ref{thm4: deg4 4*2 ram}. Still by the Weil pairing, the image of $\Gal(\bQ(E[4])/\bQ(i))$ under $\rho$ is generated by $\begin{pmatrix}
1 & 2 \\
0 & 1
\end{pmatrix}$. According to our assumption, both the inertia group and the complex conjugation factor through $\Gal(\bQ(E[4])/K)$, so the image of $\Gal(\bQ(E[4])/K)$ is generated by $\begin{pmatrix}
1 & 0 \\
0 & 3
\end{pmatrix}$, and hence both $C_2[4]$ and $C_\infty[4]$ should be either $\Spa{Q}$ or $\Spa{Q+2P}$. Since neither $\Spa{Q}$ nor $\Spa{Q+2P}$ is $G_\bQ$-invariant, the maximal odd and ramified $G_\bQ$-invariant submodule of $E[2^\infty]$ is $\Spa{Q}$, which has order $2$. Therefore, by proposition \ref{prop4: mu=m}, $\mu_E=1$.
\end{Proof}

\subsection{Computing $\mu_E$ when $E(\bQ)[4]\simeq \bZ/4\bZ$}
Now we assume $E(\bQ)[4]\simeq \bZ/4\bZ$. In this case, the image of $\rho$ has the form $\begin{pmatrix}
1 & b \\
0 & d
\end{pmatrix}$ with $0\leq b \leq 3$ and $d\in \{1,3\}$. Also, we note that either $\begin{pmatrix}
1 & 1 \\
0 & 3
\end{pmatrix}$ or $\begin{pmatrix}
1 & 1 \\
0 & 1
\end{pmatrix}$ lies inside $\im(\rho)$, as otherwise $E(\bQ)[2]$ will be $\bZ/2\bZ\times \bZ/2\bZ$. Therefore, $[\bQ(E[4]):\bQ]$ is either $4$ or $8$. Now we consider the two cases separately. \par
First we assume $[\bQ(E[4]):\bQ]=4$. By the Weil pairing, $\im(\rho)$ contains a matrix of form $\begin{pmatrix}
1 & b \\
0 & 3
\end{pmatrix}$, then $\begin{pmatrix}
1 & 1 \\
0 & 1
\end{pmatrix}$ is not in the image, hence $\im(\rho)$ is generated by $\begin{pmatrix}
1 & 1 \\
0 & 3
\end{pmatrix}$. Therefore, $\Gal(\bQ(E[4])/\bQ)\simeq \bZ/4\bZ$.
\begin{thm}\label{thm4: deg4 4}
Assume $E(\bQ)[4]\simeq \bZ/4\bZ$ and $[\bQ(E[4]):\bQ]=4$. Then $\mu_E\geq 2$.
\end{thm}

\begin{Proof}
Notice that there is an element in $I_2$ whose image under $\rho$ has determinant 3, then either $\begin{pmatrix}
1 & 1 \\
0 & 3
\end{pmatrix}$ or $\begin{pmatrix}
1 & 3 \\
0 & 3
\end{pmatrix}$ lies inside $\rho(I_2)$. Since $I_2$ acts on $C_2$ via the cyclotomic character, then $C_2[4]=\Spa{P+2Q}$. Similarly, since complex conjugations also has determinant 3, we conclude that $C_{\infty,v}[4]=\Spa{P+2Q}$. Since $\Spa{P+2Q}$ is $G_\bQ$-invariant, we can apply proposition \ref{prop4: mu=m} and conclude that $\mu_E\geq 2$.
\end{Proof}

Next we consider the case $[\bQ(E[4]):\bQ]=8$, then $\Gal(\bQ(E[4])/\bQ)\simeq D_4$. Then $\im(\rho)$ is generated by $r=\begin{pmatrix}
1 & 1 \\
0 & 1
\end{pmatrix}$ and $s=\begin{pmatrix}
1 & 0 \\
0 & 3
\end{pmatrix}$ (here we abuse the notations by identifying the matrices with Galois elements). Then the quadratic subextension $K$ corresponding to $\Spa{r^2,s}$ is the unique quadratic subextension with $E(K)[4]\simeq \bZ/4\bZ\times \bZ/2\bZ$. Write $K$ as $\bQ(\sqrt{m})$ with $m$ square-free. Since $E(\bQ)[2]\neq E[2]$, we have $\bQ(E[2])=K$, hence by lemma \ref{lem2: not Q(sqrt2)}, we conclude that $2\nmid m$.
\begin{thm}\label{thm4: deg8 4}
Assume $E(\bQ)[4]\simeq \bZ/4\bZ$ and $[\bQ(E[4]):\bQ]=8$. Let $K$ be the unique quadratic subfield of $\bQ(E[4])$ with $E(K)[4]\simeq \bZ/4\bZ\times \bZ/2\bZ$. If $K$ is imaginary and 2 ramifies in $K$, then $\mu_E\geq 2$. Otherwise, $\mu_E=0$.
\end{thm}

We need a proposition which will also be used frequently in the remaining cases.
\begin{prop}\label{prop4: E[2] nonsplit mu=0}
Assume $\bQ(E[2])$ is a quadratic field. Suppose $C_2[2]$ is not contained in $E(\bQ)$, then $\mu_E=0$
\end{prop}

\begin{Proof}
This is a re-statement of proposition 4.1 of \cite{LY26}, since if $\bQ(E[2])$ is quadratic, the nontrivial element of $\Gal(\bQ(E[2])/\bQ)$ acts on $E[2]$ by $\begin{pmatrix}
1 & 1 \\
0 & 1
\end{pmatrix}$ after we fix a basis $\{P,Q\}$. Therefore $C_2[2]$ is $G_\bQ$-invariant if and only if $C_2[2]=\Spa{P}$ if and only if $C_2[2]$ is contained in $E(\bQ)$.
\end{Proof}

Now we prove theorem \ref{thm4: deg8 4}. The key idea is still using the mod-4 representation to determine $C_2[4]$ and $C_{\infty,v}[4]$.
\begin{Proof}[theorem \ref{thm4: deg8 4}]
If $K$ is imaginary and 2 ramifies in $K$, say $K=\bQ(\sqrt{m})$, then $K'=\bQ(\sqrt{m})$ is real and 2 is unramified in $K'$. Since $K'$ correspond to $\Spa{r^2, rs}$, then either $\begin{pmatrix}
1 & 3 \\
0 & 3
\end{pmatrix}$ or $\begin{pmatrix}
1 & 1 \\
0 & 3
\end{pmatrix}$ inside $I_2$. Hence $C_2[4]=\Spa{P+2Q}$. By the same reason, there exists $v\mid \infty$ such that $C_{\infty,v}[4]=\Spa{P+2Q}$, which is $G_\bQ$-invariant, hence $C_\infty[4]=\Spa{P+2Q}$, hence by proposition \ref{prop4: mu=m}, we have $\mu_E\geq 2$. \par
Next we assume $K$ is real and 2 ramifies inside $K$. Since $K$ is real, the image of complex conjugations under $\rho$ lie inside $\Spa{r^2,s}$, hence it is either $\begin{pmatrix}
1 & 0 \\
0 & 3
\end{pmatrix}$ or $\begin{pmatrix}
1 & 2 \\
0 & 3
\end{pmatrix}$, which implies $C_{\infty,v}[2]$ is either $\Spa{2P+2Q}$ or $\Spa{2Q}$. On the other hand, since 2 in unramified in $K'$, then the images of $I_2$ under $\rho$ lie inside $\Spa{r^2,rs}$, hence it contains either $\begin{pmatrix}
1 & 1 \\
0 & 3
\end{pmatrix}$ or $\begin{pmatrix}
1 & 3 \\
0 & 3
\end{pmatrix}$. In either case, $C_2[2]=\Spa{2P}$. By proposition \ref{prop4: mu=m}, we are able to conclude that $\mu_E=0$. If $K$ is imaginary while 2 is unramified inside $K$, we just switch the roles of $K$ and $K'$, and the conclusion is still the same. \par
Now we assume $K$ is real and 2 is unramified inside $K$. This time the images of $I_2$ under $\rho$ contains either $\begin{pmatrix}
1 & 0 \\
0 & 3
\end{pmatrix}$ or $\begin{pmatrix}
1 & 2 \\
0 & 3
\end{pmatrix}$, hence $C_2[2]$ is either $\Spa{2P+2Q}$ or $\Spa{2Q}$, neither $G_\bQ$-invariant. Hence by proposition \ref{prop4: E[2] nonsplit mu=0}, $\mu_E=0$.
\end{Proof}

\subsection{Computing $\mu_E$ when $E(\bQ)[4]\simeq\bZ/2\bZ\times \bZ/2\bZ$}
Now we assume $E(\bQ)\simeq \bZ/2\bZ\times \bZ/2\bZ$. This time the image of $\rho(\Gal(\bQ(E[4])/\bQ)$ is of the form $\begin{pmatrix}
a & b \\
0 & d
\end{pmatrix}$ with $a,d\in \{1,3\}$ and $b\in\{0,2\}$. \par
We first consider the case $[\bQ(E[4]):\bQ]=4$, this implies $\Gal(\bQ(E[4])/\bQ)\simeq (\bZ/2\bZ)^2$. We first consider the possible Galois images of $\Gal(\bQ(E[4])/\bQ)$ under $\rho$. 
\begin{lemma}\label{lem4: Galois image}
Assume $E(\bQ)[4]\simeq(\bZ/2\bZ)^2$ and $[\bQ(E[4]):\bQ]=4$, then up to isomorphism, there are three possible images inside $\Aut(E[4])$, namely \[\mbox{Type A :} \Spa{\begin{pmatrix}
3 & 0 \\
0 & 1
\end{pmatrix},\begin{pmatrix}
3 & 2 \\
0 & 1
\end{pmatrix}},\] \[\mbox{Type B :}\Spa {\begin{pmatrix}
3 & 0 \\
0 & 1
\end{pmatrix}, \begin{pmatrix}
1 & 0 \\
0 & 3
\end{pmatrix}},\] and \[\mbox{Type C :}\Spa {\begin{pmatrix}
3 & 0 \\
0 & 1
\end{pmatrix},\begin{pmatrix}
1 & 2 \\
0 & 3
\end{pmatrix}}.\]
\end{lemma}

\begin{Proof}
For simplicity, set $a=\begin{pmatrix}
3 & 0 \\
0 & 1
\end{pmatrix}, b=\begin{pmatrix}
1 & 2 \\
0 & 1
\end{pmatrix}$ and $c=\begin{pmatrix}
1 & 0 \\
0 & 3
\end{pmatrix}$, then there are 7 subgroups of order 4 inside $\Spa{a,b,c}$. However, we need to exclude the case $\Spa{b,c}$, as it satisfies $E(\bQ)[4]\simeq \bZ/4\bZ\times \bZ/2\bZ$, which contradicts to our assumption. Also, by the Weil pairing, there exists a matrix with determinant 3, hence it cannot be $\Spa{b,ac}$.\par
In the remaining cases, $\Spa{a,b}$ is the type A above, $\Spa{a,c}$ is the type B and $\Spa{a,bc}$ is the type C. Considering the change of basis $P\mapsto P, Q\mapsto P+Q$, we can see that $\Spa{c,ab}$ is isomorphic to type C while $\Spa{ac,bc}$ is isomorphic to type B.
\end{Proof}

Now we deal with the three cases separately. For type A, complex conjugations act on $E[4]$ by either $\begin{pmatrix}
3 & 0 \\
0 & 1
\end{pmatrix}$ or $\begin{pmatrix}
3 & 2 \\
0 & 1
\end{pmatrix}$, so $C_{\infty,v}[4]$ is either $\Spa{P}$ or $\Spa{P+2Q}$, both $G_\bQ$-invariant, hence $C_{\infty,v}$ is independent of the infinite prime $v$ and we denote it by $C_\infty$.

\begin{thm}\label{thm4: deg4 2*2 Type A}
Suppose $E(\bQ)[4]\simeq \bZ/2\bZ\times \bZ/2\bZ$ and $[\bQ(E[4]):\bQ]=4$. Further assume the Galois image under $\rho$ is of type A. If the intersection of $C_2[4]$ and $C_\infty[4]$ has order 2, then $\mu_E=1$. If $C_2[4]=C_\infty[4]$, then $\mu_E\geq 2$.
\end{thm}

We first prove a proposition which will also be useful in the remaining cases.
\begin{prop}\label{prop4: C2 intersect Cinfty}
Assume $C_2[4]$ and $C_\infty[4]$ are both $G_\bQ$-invariant. If $C_2[4]=C_\infty[4]$, then $\mu_E\geq 2$. If the intersection of $C_2[4]$ and $C_\infty[4]$ has order 2, then $\mu_E=1$.
\end{prop}

\begin{Proof}
According to our assumptions, both parts follow directly from proposition \ref{prop4: mu=m}.
\end{Proof}

\begin{Proof}[theorem \ref{thm4: deg4 2*2 Type A}]
In order to apply proposition \ref{prop4: C2 intersect Cinfty}, we just need to show that both $C_2[4]$ and $C_{\infty,v}[4]$ are $G_\bQ$-invariant. We have already explained before the statement of theorem \ref{thm4: deg4 2*2 Type A} that $C_\infty[4]$ is $G_\bQ$-invariant and $C_\infty[4]$ is either $\Spa{P}$ or $\Spa{P+2Q}$. Similarly, since $I_2$ contains either $\begin{pmatrix}
3 & 0 \\
0 & 1
\end{pmatrix}$ or $\begin{pmatrix}
3 & 2 \\
0 & 1
\end{pmatrix}$, we can see that $C_2[4]$ is either $\Spa{P}$ or $\Spa{P+2Q}$, which is also $\G_\bQ$-invariant. Now the intersection of $C_2[4]$ and $C_\infty[4]$ is nontrivial, hence the theorem follows from proposition \ref{prop4: C2 intersect Cinfty}.
\end{Proof}

Next we consider type B. This time complex conjugations act on $E[4]$ by either $\begin{pmatrix}
3 & 0 \\
0 & 1
\end{pmatrix}$ or $\begin{pmatrix}
1 & 0 \\
0 & 3
\end{pmatrix}$, hence $C_{\infty,v}[4]$ is one of $\Spa{P}$, $\Spa{P+2Q}$, $\Spa{Q}$ or $\Spa{Q+2P}$. All of them are $G_\bQ$-invariant, so we can still denote them by $C_\infty[4]$.

\begin{thm}\label{thm4: deg4 2*2 Type B}
Suppose $E(\bQ)[4]\simeq \bZ/2\bZ\times \bZ/2\bZ$ and $[\bQ(E[4]):\bQ]=4$. Further assume the Galois image under $\rho$ is of type B. Let $K$ be the unique real quadratic subfield of $\bQ(E[4])$.Then\begin{enumerate} 
\item If 2 ramifies inside $K$, then $\mu_E=0$. 
\item If 2 is unramified in $K$ and the intersection of $C_2[4]$ and $C_\infty[4]$ has order 2, then $\mu_E=1$.
\item If 2 is unramified in $K$ and $C_2[4]=C_\infty[4]$, then $\mu_E\geq 2$. 
\end{enumerate}
\end{thm}

\begin{Proof}
By the Weil pairing, the nontrivial element of $\Gal(\bQ(E[4])/\bQ(i))$ acts on $E[4]$ by $\begin{pmatrix}
3 & 0 \\
0 & 3
\end{pmatrix}$, hence by switching $P$ and $Q$ if necessary, we may assume the nontrivial element of $\Gal(\bQ(E[4])/K)$ acts by $\begin{pmatrix}
3 & 0 \\
0 & 1
\end{pmatrix}$, hence complex conjugations act on $E[4]$ by$\begin{pmatrix}
3 & 0 \\
0 & 1
\end{pmatrix}$, therefore $C_\infty[4]$ is either $\Spa{P}$ or $\Spa{P+2Q}$. \par
Say $K=\bQ(\sqrt{m})$, let $K'=\bQ(\sqrt{-m})$. If 2 ramifies in $K$, then it is unramified in $K'$ as the ramification index is $e=2$. Since the image of $\Gal(\bQ(E[4])/K')$ under $\rho$ is generated by $\begin{pmatrix}
1 & 0 \\
0 & 3
\end{pmatrix}$, then $I_2$ acts on $E[4]$ by $\begin{pmatrix}
1 & 0 \\
0 & 3
\end{pmatrix}$, then $C_2[4]$ is either $\Spa{Q}$ or $\Spa{Q+2P}$. Now $C_2[4]$ and $C_\infty[4]$ have trivial intersection, so by proposition \ref{prop4: mu=m}, we conclude that $\mu_E=0$. \par
If 2 is unramified in $K$, then $I_2$ acts on $E[4]$ via $\begin{pmatrix}
3 & 0 \\
0 & 1
\end{pmatrix}$, hence $C_2[4]$ is either $\Spa{P}$ or $\Spa{P+2Q}$. The remaining part of the theorem follows from proposition \ref{prop4: C2 intersect Cinfty}.
\end{Proof}

Finally we consider type C. For simplicity, let $r=\begin{pmatrix}
3 & 0 \\
0 & 1
\end{pmatrix}$ and $s=\begin{pmatrix}
1 & 2 \\
0 & 3
\end{pmatrix}$. Then $\Gal(\bQ(E[4])/\bQ(i))$ correspond to $\Spa{rs}$. Let $K$ be the quadratic subfield corresponding to the subgroup $\Spa{r}$, which can be characterized as the unique quadratic subfield with a 4-torsion that is not $G_\bQ$-invariant. Say $K=\bQ(\sqrt{m})$, set $K'=\bQ(\sqrt{-m})$. Using an argument close to lemma \ref{lem4: e=2}, we know that $2\nmid m$, hence 2 is unramified in exactly one of $K$ and $K'$.

\begin{thm}\label{thm4: deg4 2*2 Type C}
Suppose $E(\bQ)[4]\simeq \bZ/2\bZ\times \bZ/2\bZ$ and $[\bQ(E[4]):\bQ]=4$. Further assume the Galois image under $\rho$ is of type C. Let $K$ be the unique quadratic subfield with a 4-torsion that is not $G_\bQ$-invariant. Then \begin{enumerate}
\item If 2 ramifies inside $K$ with $K$ real, then $\mu_E=0$.
\item If 2 is unramified in $K$ with K imaginary, then $\mu_E=0$.
\item If 2 is ramified in $K$ with $K$ imaginary, then $\mu_E=1$. 
\item If 2 is unramified in $K$ with $K$ real, and the intersection of $C_2[4]$ and $C_\infty[4]$ has order 2 then $\mu_E=1$.
\item If 2 is unramified in $K$ with $K$ real and $C_2[4]=C_\infty[4]$, then $\mu_E\geq 2$. 
\end{enumerate}
\end{thm}

The main difficulty lies inside the case where 2 is ramified in $K$ with $K$ imaginary. We need a proposition for this case, which can be viewed as a generalization of the proof of theorem \ref{thm4: deg4 4*2 unram}.

\begin{prop}\label{prop4: E[4] nonsplit mu<=1}
Assume $E(\bQ)[4]\simeq (\bZ/2\bZ)^2$. Let $P,Q$ be a basis of $E[4]$ with $\Spa{P}$ a $G_\bQ$-invariant submodule. Suppose $C_2[4]=\Spa{Q}$ is not $G_\bQ$-invariant. Then $\mu_E\leq 1$.
\end{prop}

\begin{Proof}
It follows directly from proposition \ref{prop4: mu=m}. As $C_2[4]$ is not $G_\bQ$-invariant, the maximal odd and ramified $G_\bQ$-invariant submodule of $E[2^\infty]$ has order $\leq 2$.
\end{Proof}

Now we are able to prove theorem \ref{thm4: deg4 2*2 Type C}.
\begin{Proof}[theorem \ref{thm4: deg4 2*2 Type C}]
We first consider the case where 2 is ramfied in $K$ while $K$ is real. Now complex conjugation acts by $\begin{pmatrix}
3 & 0 \\
0 & 1
\end{pmatrix}$, and hence $C_{\infty,v}[4]$ is either $\Spa{P}$ or $\Spa{P+2Q}$, both $G_\bQ$-invariant. In either case, $C_\infty[2]=\Spa{2P}$. On the other hand, since 2 is unramified in $K'$, the image of $I_2$ is generated by $\begin{pmatrix}
1 & 2 \\
0 & 3
\end{pmatrix}$, and $C_2[4]$ is either $\Spa{P+Q}$ or $\Spa{P+3Q}$. In either case, $C_2[2]=\Spa{2P+2Q}$, hence by proposition \ref{prop4: mu=m}, $\mu_E=0$. If 2 is unramified in $K$ while $K$ is imaginary, we just swap the roles of $K$ and $K'$.\par
Next assume 2 is unramified in $K$ and $K$ is real. Then both complex conjugations and $I_2$ act via $\begin{pmatrix}
3 & 0 \\
0 & 1
\end{pmatrix}$, which means both $C_2[4]$ and $C_\infty[4]$ are $G_\bQ$-invariant with nontrivial intersection, hence the result follows from proposition \ref{prop4: C2 intersect Cinfty}.\par
Finally assume 2 is ramified in $K$ and $K$ is imaginary, then 2 is unramified in $K'$ with $K'$ real. This time $C_2[2]=\Spa{2P+2Q}=C_{\infty,v}[2]$ which is $G_\bQ$-invariant while $C_2[4]$ is not $G_\bQ$-invariant. Now we apply a change of basis by letting $P$ and a generator of $C_2[4]$, say $R$, to be a basis of $E[4]$, then by direct computation we notice that $G_{K'}$ acts trivially on $\Spa{2P,R}/\Spa{2R}$, hence by proposition \ref{prop4: E[4] nonsplit mu<=1}, we have $\mu_E\leq 1$. On the other hand, since both $C_2[2]$ and $C_\infty[2]$ are $\Spa{2P+2Q}$, by proposition \ref{prop4: mu=m}, we have $\mu_E\geq 1$. Combining the two results, we conclude that $\mu_E=1$.
\end{Proof}

Next we consider the case $[\bQ(E[4]):\bQ]=8$, this time $\rho(\Gal(\bQ(E[4])/\bQ)$ is generated by $r=\begin{pmatrix}
3 & 0 \\
0 & 1
\end{pmatrix}$, $s=\begin{pmatrix}
1 & 2 \\
0 & 1
\end{pmatrix}$ and $t=\begin{pmatrix}
1 & 0 \\
0 & 3
\end{pmatrix}$, hence $\Gal(\bQ(E[4])/\bQ)\simeq (\bZ/2\bZ)^3$. Let $K$ be the quadratic subfield of $\bQ(E[4])$ corresponding to $\Spa{s,t}$, the subfield characterized as the unique quadratic subfield with $E(K)[4]\simeq \bZ/4\bZ\times \bZ/2\bZ$. Say $K=\bQ(\sqrt{m})$. By the Weil pairing, $\bQ(i)$ corresponds to $\Spa{rt,s}$, hence $\bQ(i,\sqrt{m})$ corresponds to $\Spa{s}$ and $\bQ(\sqrt{-m})$ corresponds to $\Spa{r,s}$. Since $[2]$ on the formal group of $E$ is a quadratic polynomial, then $\bQ(E[4])/\bQ$ has ramification index $e\leq 2$, hence $2\nmid m$. \par 
Notice that the two biquadratic subfields of $E[4]$ that correspond to  $\Spa{r}$ and $\Spa{rs}$ are the two fields characterized as the only biquadratic subfields in $\bQ(E[4])$ with no subfield with 4-torsion isomorphic to $\bZ/4\bZ\times \bZ/2\bZ$. Since both of these two fields are above $\bQ(\sqrt{-m})$ and neither of them contains $\bQ(i)$, so they can be written as $\bQ(\sqrt{-m}, \sqrt{n})$ and $\bQ(\sqrt{-m},\sqrt{-n})$ with $n$ a positive integer not divisible by $m$.
\begin{thm}\label{thm4: deg8 2*2}
Suppose $E(\bQ)[4]\simeq \bZ/2\bZ\times \bZ/2\bZ$ and $[\bQ(E[4]):\bQ]=8$. Let $K=\bQ(\sqrt{m})$ be the unique quadratic subfield with $E(K)[4]\simeq \bZ/4\bZ\times \bZ/2\bZ$, and $L=\bQ(\sqrt{-m},\sqrt{n})$ be one of the two biquadratic subfields (the other is $\bQ(\sqrt{-m},\sqrt{-n})$) with \begin{enumerate}
\item $E(L)[4]\simeq \bZ/4\bZ\times \bZ/2\bZ$.
\item For any subfield $F\subset L$, $E(F)[4]\subsetneq E(L)[4]$.
\item $L$ contains a real quadratic subfield $L'=\bQ(\sqrt{n})$ with $m\nmid n$.
\end{enumerate}  Then
\begin{enumerate}
\item If $K$ is real and 2 ramifies in $K$, then $\mu_E=0$.
\item If $K$ is imaginary and 2 is unramfied in $K$, then $\mu_E=0$.
\item If $K$ is real and 2 is unramified in $K$ while ramified in $L'$, then $\mu_E=0$.
\item If $K$ is real and 2 is unramified in both $K$ and $L'$, then $\mu_E=1$.
\item If $K$ is imaginary and 2 ramifies in $K$ with the intersection of $C_2[4]$ and $C_\infty[4]$ has order 2 , then $\mu_E=1$.
\item If $K$ is imaginary and 2 ramifies in $K$ with $C_2[4]=C_\infty[4]$ has order 4 , then $\mu_E\geq 2$.
\end{enumerate}
\end{thm}

\begin{Proof}
If $K$ is real and 2 ramifies in $K$, then the complex conjugations act on $E[4]$ by either $\begin{pmatrix}
1 & 0 \\
0 & 3
\end{pmatrix}$ or $\begin{pmatrix}
1 & 2 \\
0 & 3
\end{pmatrix}$, hence $C_{\infty,v}[2]$ is either $\Spa{2Q}$ or $\Spa{2P+2Q}$. On the other hand, since 2 is unramified in $K'$, $I_2$ acts on $E[4]$ by either $\begin{pmatrix}
3 & 0 \\
0 & 1
\end{pmatrix}$ or $\begin{pmatrix}
3 & 2 \\
0 & 1
\end{pmatrix}$, hence $C_2[2]=\Spa{2P}$. By proposition \ref{prop4: mu=m}, we conclude that $\mu_E=0$. If $K$ is imaginary and 2 is unramified in $K$, we just switch the roles of $K$ and $K'$. \par
If both $K$ and $L'$ are real and 2 is unramified in $K$ while ramified in $L'$, we first show that it doesn't depend on the choice of $L$. Notice the assumptions imply that 2 is unramified in $\bQ(\sqrt{m})$ while ramified in $\bQ(\sqrt{n})$, therefore 2 is ramified in $\bQ(\sqrt{mn})$, hence it is independent of the choice if we replace $n$ by $-n$. Now complex conjugations act on $E[4]$ via the image of the nontrivial element of $\Gal(\bQ(E[4])/\bQ(\sqrt{m},\sqrt{n})$. On the other hand, 2 is unramified in $\bQ(\sqrt{m},\sqrt{-n})$, hence the image of $I_2$ under $\rho$ is exactly the image of $\Gal(\bQ(E[4])/\bQ(\sqrt{m},\sqrt{-n}))$ under $\rho$. Since $K$ corresponds to $\Spa{s,t}$, one of them is $\Spa{t}$ and the other is $\Spa{st}$. Hence one of them contains $\begin{pmatrix}
1 & 0 \\
0 & 3
\end{pmatrix}$ while the other one contains $r=\begin{pmatrix}
1 & 2 \\
0 & 3
\end{pmatrix}$, hence one of $C_2[2]$ and $C_{\infty,v}[2]$ is $\Spa{2P+2Q}$ while the other one is $\Spa{2Q}$. Now the result follows from proposition \ref{prop4: mu=m}. \par
If $K$ is real and 2 is unramified in both $K$ and $L'$, we first show that it doesn't depend on the choice of $L$. Notice the assumptions imply that 2 is unramified in both $\bQ(\sqrt{m})$ and $\bQ(\sqrt{n})$, therefore 2 is unramified in $\bQ(\sqrt{mn})$, hence it is independent of the choice if we replace $n$ by $-n$. Now both $C_2[2]$ and $C_{\infty,v}[2]$ are equal, which is either $\Spa{2P+2Q}$ or $\Spa{2Q}$, hence by proposition \ref{prop4: mu=m}, we have $\mu_E\geq 1$. On the other hand, notice that $C_2[4]$ is one of $\Spa{P+Q}$, $\Spa{P+3Q}$, $\Spa{Q}$ or $\Spa{Q+2P}$, none of which is $G_\bQ$-invariant. Therefore, by proposition \ref{prop4: E[4] nonsplit mu<=1}, we have $\mu_E\leq 1$. Combining both sides, we conclude that $\mu_E=1$.\par
If $K$ is imaginary and 2 ramifies in $K$, then $K'$ is real and 2 is unramified in $K'$, hence complex conjugations and elements in $I_2$ act on $E[4]$ via either $\begin{pmatrix}
3 & 0 \\
0 & 1
\end{pmatrix}$ or $\begin{pmatrix}
3 & 2 \\
0 & 1
\end{pmatrix}$. In either case, $C_2[2]=C_{\infty,v}[2]=\Spa{2P}$. Hence $C_{\infty,v}$ is indepedent of $v$. Now since both $\Spa{P}$ and $\Spa{P+2Q}$ are $G_\bQ$-invariant, we may apply proposition \ref{prop4: C2 intersect Cinfty}.
\end{Proof}

\subsection{Computing $\mu_E$ when $E(\bQ)[4]\simeq\bZ/2\bZ$}
Now we consider the case $E(\bQ)[4]\simeq \bZ/2\bZ$. This time the image of $\Gal(\bQ(E[4])/\bQ)$ under $\rho$ is of the form $\begin{pmatrix}
a & b \\
0 & d
\end{pmatrix}$ with $a\in\{1,3\}$, $b\in \{0,1,2,3\}$ and $d\in \{1,3\}$. Hence $[\bQ(E[4]):\bQ]=8$ or 16. Also, locally we can see that $[\bQ_2(E[4]):\bQ_2]=8$ with inertia degree either $e=2$ or $e=4$, depending on whether the 2-torsion in $\bQ$ is in $C_2$ or not. \par
We first assume $[\bQ(E[4]):\bQ]=8$. By changing a basis if necessary, we may assume $\im(\rho)$ contains $r=\begin{pmatrix}
1 & 1 \\
0 & 1
\end{pmatrix}$. Since $E(\bQ)[4]\simeq \bZ/2\bZ$, $\im(\rho)$ also contains $s=\begin{pmatrix}
3 & 0 \\
0 & 1
\end{pmatrix}$, and $\Gal(\bQ(E[4])/\bQ)\simeq D_4$. \par
Now $\Gal(\bQ(E[4])/\bQ)$ contains 3 quadratic subfields, namely the subextensions corresponding to $\Spa{s,r^2}$, $\Spa{r}$ and $\Spa{sr,r^2}$. By the Weil pairing, the one corresponding to $\Spa{r}$ is $\bQ(i)$. Let $K=\bQ(\sqrt{m})$ be the one corresponding to $\Spa{s,r^2}$, which is the unique quadratic subfield with $E(K)[4]\simeq \bZ/2\bZ\times \bZ/2\bZ$. We remark that $2\nmid m$ : Suppose 2 has ramification index $e=2$, then $2\nmid m$ as otherwise, 2 is totally ramified in $\bQ(\sqrt{m},i)\subset \bQ(E[4])$. Now suppose 2 has ramification index $e=4$, then there exists a quadratic subfield $F$ of $\bQ(E[4])$, either $\bQ(\sqrt{m})$ or $\bQ(\sqrt{-m})$ such that 2 is unramified in $F$. This still implies $2\nmid m$. Therefore, we conclude that 2 is unramified in one of $K=\bQ(\sqrt{m})$ and $K'=\bQ(\sqrt{-m})$.
\begin{thm}\label{thm4: deg8 2}
Suppose $E(\bQ)[4]\simeq \bZ/2\bZ$ and $[\bQ(E[4]):\bQ]=8$. Let $K=\bQ(\sqrt{m})$ be the unique quadratic subfield with $E(K)\simeq \bZ/2\bZ\times \bZ/2\bZ$. Then 
\begin{enumerate}
\item If either $K$ is real or 2 is unramified in $K$, and the intersection of $C_2[4]$ and $C_\infty[4]$ has order 2, then $\mu_E=1$.
\item If either $K$ is real or 2 is unramified in $K$ and  $C_2[4]=C_\infty[4]$, then $\mu_E\geq2$.
\item If $K$ is imaginary and 2 is ramified in $K$, then $\mu_E\geq 2$.
\end{enumerate}
\end{thm}

\begin{Proof}
We first make some general observations. If $K$ is real, then complex conjugations act on $E[4]$ by either $\begin{pmatrix}
3 & 0 \\
0 & 1
\end{pmatrix}$ or $\begin{pmatrix}
3 & 2 \\
0 & 1
\end{pmatrix}$, hence $C_{\infty,v}[4]$ is either $\Spa{P}$ or $\Spa{P+2Q}$, both of which are $G_\bQ$-invariant. Therefore we may denote them by $C_\infty[4]$. On the other hand, if $K$ is imaginary, then $K'=\bQ(\sqrt{-m})$ is real, hence complex conjugations act on $E[4]$ by either $\begin{pmatrix}
3 & 1 \\
0 & 1
\end{pmatrix}$ or $\begin{pmatrix}
3 & 3 \\
0 & 1
\end{pmatrix}$, hence $C_{\infty,v}[4]=\Spa{P}$, which is still $G_\bQ$-invariant. Similarly, by considering the possible images of $\rho(I_2)$, we have $C_2[4]=\Spa{P+2Q}$ or $\Spa{P}$ if 2 is unramified in $K$, while $C_2[4]=\Spa{P}$ if 2 is ramified in $K$. In either case, $C_2[4]$ is $G_\bQ$-invariant.\par
Now if either $K$ is real or 2 is unramified in $K$, then both $C_2[4]$ and $C_\infty[4]$ can be $\Spa{P}$ or $\Spa{P+2Q}$, so we get the conclusion from proposition \ref{prop4: C2 intersect Cinfty}. If $K$ is imaginary and 2 is ramified in $K$, then $C_2[4]=C_\infty[4]=\Spa{P}$, hence by proposition \ref{prop4: mu=m}, $\mu_E\geq 2$.
\end{Proof}

Finally we consider the case where $[\bQ(E[4]):\bQ]=16$, and this time the image of $\Gal(\bQ(E[4])/\bQ)$ under $\rho$ is generated by $r=\begin{pmatrix}
1 & 1 \\
0 & 1
\end{pmatrix}$, $s=\begin{pmatrix}
1 & 0 \\
0 & 3
\end{pmatrix}$ and $\begin{pmatrix}
3 & 0 \\
0 & 1
\end{pmatrix}$. By the Weil pairing, $\bQ(i)$ corresponds to the subgroup $\Spa{r,ts}$. Now let $K=\bQ(\sqrt{m})$ be the quadratic subfield of $\bQ(E[4])$ corresponding to $\Spa{s,r^2,t}$, which is the unique quadratic subfield with $E(K)[4]\simeq (\bZ/2\bZ)^2$. Let $F=\bQ(\sqrt{m},\sqrt{n})$ be the biquadratic subfield corresponding to $\Spa{r^2,s}$. This is the unique biquadratic subfield with $E(F)[4]\simeq \bZ/4\bZ\times \bZ/2\bZ$. Note that $F$ has two quadratic subfields corresponding to $\Spa{r,s}$ and $\Spa{s,r^2,tr}$, which are the only quadratic subfields with 4-torsions. We denote the one corresponding to $\Spa{r,s}$ by $K=\bQ(\sqrt{m})$, and the other by $K'=\bQ(\sqrt{n})$. Let $L$ be the quadratic subfield corresponding to $\Spa{s,r^2,t}$, the unique quadratic subfield with $E(L)[4]\simeq (\bZ/2\bZ)^2$. From direct computation, we can see that $L=\bQ(\sqrt{mn})$. Now using Galois theory, we can see that $F_1:=\bQ(\sqrt{m},\sqrt{-n})$ corresponds to $\Spa{r^2,rs}$, $F_2:=\bQ(\sqrt{-m},\sqrt{n})$ corresponds to $\Spa{r^2,tr}$, and $F_3:=\bQ(\sqrt{-m},\sqrt{-n})$ corresponds to $\Spa{r^2,t}$. \par
Before stating the theorem, we remark that neither $m$ nor $n$ is divisible by 2: Suppose instead $2\mid m$, then 2 is totally ramified in $\bQ(\sqrt{m},i)$, hence $[\bQ_2(\sqrt{m},i):\bQ_2]=4$. However, since $[\bQ_2(E[4]):\bQ_2]\leq 8$, which means 2 should split in $\bQ(\sqrt{mn})$ as this is the only quadratic field with the full 2-torsion, contradiction arises. Similarly, we have $2\nmid n$. This implies there is exactly one field among $F$, $F_1$, $F_2$ and $F_3$ such that 2 is unramified there.
\begin{thm}\label{thm4: deg16 2}
Suppose $E(\bQ)[4]\simeq \bZ/2\bZ$ and $[\bQ(E[4]):\bQ]=16$. Let $F$ be the unique biquadratic subfield with $E(F)[4]\simeq \bZ/4\bZ\times \bZ/2\bZ$ and $L$ be the unique quadratic subfield with $E(L)[4]\simeq \bZ/2\bZ\times \bZ/2\bZ$. Then 
\begin{enumerate}
\item If $F$ is real or 2 is unramified in $F$, then $\mu_E=0$.
\item If $F$ is imaginary and 2 is ramified in $F$ while $L$ is real or unramified, further assume the intersection of $C_2[4]$ and $C_\infty[4]$ has order 2, then $\mu_E=1$. 
\item If $F$ is imaginary and 2 is ramified in $F$ while $L$ is real or unramified, further assume $C_2[4]=C_\infty[4]$, then $\mu_E\geq 2$.
\item If $L$ is imaginary and 2 is ramified in $L$. Let $F'$ be the unique real quadratic subfield of $F$. If 2 is ramified in $F'$, then $\mu_E=1$.
\item If $L$ is imaginary and 2 is ramified in $L$. Let $F'$ be the unique real quadratic subfield of $F$. If 2 is unramified in $F'$, then $\mu_E\geq 2$.
\end{enumerate}
\end{thm}

\begin{Proof}
The key ingredient of the proof is still analyzing the possibilities for $C_2[4]$ and $C_\infty[4]$. We first consider the case where $F$ is real. This time complex conjugations act on $E[4]$ by $\begin{pmatrix}
1 & 0 \\
0 & 3
\end{pmatrix}$ or $\begin{pmatrix}
1 & 2 \\
0 & 3
\end{pmatrix}$, hence $C_{\infty,v}[2]$ is either $\Spa{2P+2Q}$ or $\Spa{2Q}$. If $C_2[2]=\Spa{2P}$, then by proposition \ref{prop4: mu=m}, we have $\mu_E=0$. If $C_2[2]=\Spa{2P+2Q}$ or $\Spa{2Q}$, since neither of them is $G_\bQ$-invariant, we conclude from proposition \ref{prop4: E[2] nonsplit mu=0} that $\mu_E=0$. \par
Before we move on, we first make some general observations when one of $F_1$ to $F_3$ is totally real. If $F_1$ is totally real, then complex conjugations act on $E[4]$ via $\begin{pmatrix}
1 & 1 \\
0 & 3
\end{pmatrix}$ or $\begin{pmatrix}
1 & 3 \\
0 & 3
\end{pmatrix}$, hence $C_{\infty,v}[4]=\Spa{P+2Q}$, which is $G_\bQ$-invariant, and we can denote it by $C_\infty[4]$. If $F_2$ is totally real, then the image of complex conjugations are $\begin{pmatrix}
1 & 1 \\
0 & 3
\end{pmatrix}$ or $\begin{pmatrix}
1 & 3 \\
0 & 3
\end{pmatrix}$, hence $C_\infty[4]=\Spa{P}$. If $F_3$ is totally real, then the images of complex conjugations are $\begin{pmatrix}
3 & 0 \\
0 & 1
\end{pmatrix}$ and $\begin{pmatrix}
3 & 2 \\
0 & 1
\end{pmatrix}$, hence $C_\infty[4]$ can be either $\Spa{P}$ or $\Spa{P+2Q}$. Similarly we can do the same study towards the images of $I_2$. We have $C_2[4]=\Spa{P}$ if 2 is unramified in $F_1$, $C_2[4]=\Spa{P+2Q}$ if 2 is unramified in $F_2$ and $C_2[4]$ can be either $\Spa{P}$ or $\Spa{P+2Q}$ if 2 is unramified in $F_3$. \par
Now we consider the case where 2 is unramified in $F$ and $F$ is imaginary. Using a similar argument, we can show that $C_2[2]=\Spa{2P+2Q}$ or $\Spa{2Q}$. However, since one of $F_1$ to $F_3$ is totally real, $C_\infty[2]=\Spa{P}$ or $\Spa{P+2Q}$. Hence still by proposition \ref{prop4: mu=m}, we have $\mu_E=0$. \par
Next we consider the case where $F$ is not real while $L$ is real. This implies $F_3$ is totally real, hence $C_{\infty,v}[4]$ is either $\Spa{P}$, or $\Spa{P+2Q}$. Now since 2 is ramified in $F$, then it is unramified in one of $F_1$ to $F_3$, hence $C_2[4]$ can be either $\Spa{P}$ or $\Spa{P+2Q}$. Now we are in the situation where proposition \ref{prop4: C2 intersect Cinfty} is applicable. Similarly, if 2 is unramified in $L$, then it is unramified in $F_3$, hence $C_2[4]$ can be either $\Spa{P}$ or $\Spa{P+2Q}$. Since $F$ is imaginary, one of $F_1$ to $F_3$ is totally real, hence $C_\infty[4]$ is either $\Spa{P}$ or $\Spa{P+2Q}$. Therefore, we still can apply proposition \ref{prop4: C2 intersect Cinfty}. \par
Next we assume $L$ is imaginary and ramified while 2 is ramified in $F'$. This implies 2 is unramified in one of $F_1$ or $F_2$, while the other one is totally real. If 2 is unramified in $F_1$ while $F_2$ is totally real, then $C_2[4]=\Spa{P}$ and $C_\infty[4]=\Spa{P+2Q}$, which means their intersection has order 2, hence we can still apply proposition \ref{prop4: C2 intersect Cinfty}. The other case can be proved similarly. \par
Finally we assume $L$ is imaginary and ramified while 2 is unramified in $F'$. This implies 2 is unramfied in one of $F_1$ or $F_2$, and that field is totally real. If 2 is unramified in $F_1$ and $F_1$ is totally real, then $C_2[4]=C_\infty[4]=\Spa{P}$, hence we have $\mu_E\geq 2$ by proposition \ref{prop4: C2 intersect Cinfty}. The other case is similar.
\end{Proof}
\section{Examples}\label{sec 5}
\large
In this section, we give some examples and some directions for further developments.
\subsection{Examples with $\mu_E=0$}
We first offer some examples for elliptic curves $E$ defined over $\bQ$ with $E[4]$ reducible as a $G_\bQ$-module and  $\mu_E=0$ that can be computed using the results in our work.

\begin{exam}
Let $E: y^2+xy=x^3-4x-1$. This is 21A5 in the LMFDB table. We have $E(\bQ)[4]\simeq \bZ/4\bZ\times \bZ/2\bZ$ and $\bQ(E[4])=\bQ(i,\sqrt{7})$, hence the unique real quadratic subfield is $\bQ(\sqrt{7})$. Since 2 is ramified in $\bQ(\sqrt{7})$, by theorem \ref{thm4: deg4 4*2 ram}, we conclude that $\mu_E=0$.
\end{exam}

\begin{exam}
Let $E: y^2+xy=x^3+x^2-146x+621$. This is 33A1 in the LMFDB table. We have $E(\bQ)[4]\simeq \bZ/4\bZ$ and $[\bQ(E[4]):\bQ]=8$. Let $K=\bQ(\sqrt{3})$, then we have $\bQ(K)[4]\simeq \bZ/4\bZ\times \bZ/2\bZ$. Since $K$ is real, then according to theorem \ref{thm4: deg8 4}, we conclude that $\mu_E=0$.
\end{exam}

\begin{exam}
Let $E:y^2+xy=x^3-x^2-1773x-5720$. This is 441F5 in the LMFDB table. We have $E(\bQ)[4]\simeq \bZ/2\bZ$ and $[\bQ(E[4]):\bQ]=8$. Now $K=\bQ(\sqrt{21})$ is the unique quadratic subfield with $E(K)[4]\simeq \bZ/4\bZ\times \bZ/2bZ$, $L=\bQ(\sqrt{-3},\sqrt{7})$ is one of the two biquadratic subfields of $\bQ(E[4])$ with $E(L)[4]\simeq \bZ/4\bZ\times \bZ/2\bZ$ and lies above $\bQ(\sqrt{-21})$. Now $L'=\bQ(\sqrt{7})$. Since 2 is unramified in $K$ and $K$ is real, while 2 is ramified in $L'$. Using theorem \ref{thm4: deg8 2*2}, we conclude that $\mu_E=0$.
\end{exam}

\begin{exam}
Let $E: y^2+xy=x^3-x^2-5$. This is 45A7 in the LMFDB table. We have $E(\bQ)[4]\simeq \bZ/2\bZ$ and $[\bQ(E[4]):\bQ]=16$. Now $F=\bQ(\sqrt{-3},\sqrt{5})$ is the unique biquadratic subfield of $\bQ(E[4])$ such that $E(F)[4]\simeq \bZ/4\bZ\times \bZ/2\bZ$. Since 2 is unramified in $F$, by theorem \ref{thm4: deg16 2}, we have $\mu_E=0$.
\end{exam}

\subsection{Examples with $\mu_E=1$}
Next we offer some examples for elliptic curves $E$ defined over $\bQ$ with $E[4]$ reducible as a $G_\bQ$-module and $\mu_E=1$.

\begin{exam}
Let $E: y^2+xy+y=x^3+x^2-39x+36$. This is 231A3 in the LMFDB table. We have $E(\bQ)[4]\simeq \bZ/4\bZ\times \bZ/2\bZ$ and $\bQ(E[4])=\bQ(i,\sqrt{33})$, hence the unique real quadratic subfield is $\bQ(\sqrt{33})$. Since 2 is ramified in $\bQ(\sqrt{33})$, by theorem \ref{thm4: deg4 4*2 unram}, we conclude that $\mu_E=1$.
\end{exam}


\begin{exam}
Let $E: y^2+xy=x^3+x^2-2875x+49000$. This is 975I5 in the LMFDB table. We have $E(\bQ)[4]\simeq \bZ/2\bZ$ and $[\bQ(E[4]):\bQ]=8$. Now $K=\bQ(\sqrt{5})$ is the unique quadratic subfield with $E(K)[4]\simeq \bZ/4\bZ\times \bZ/2bZ$, $L=\bQ(\sqrt{-5},\sqrt{13})$ is one of the two biquadratic subfields of $\bQ(E[4])$ with $E(L)[4]\simeq \bZ/4\bZ\times \bZ/2\bZ$ and lies above $\bQ(\sqrt{-5})$. Now $L'=\bQ(\sqrt{13})$. Since 2 is unramified in both $K$ and $L'$ with $K$ real, by theorem \ref{thm4: deg8 2*2}, we conclude that $\mu_E=1$.
\end{exam}

\begin{exam}
Let $E:y^2+xy=x^3-x^2-990x+22765$. This is 45A3 in the LMFDB table. We have $E(\bQ)[4]\simeq \bZ/2\bZ$ and $[\bQ(E[4]):\bQ]=16$. Now $F=\bQ(\sqrt{3},\sqrt{-5})$ is the unique biquadratic subfield of $\bQ(E[4])$ such that $E(F)[4]\simeq \bZ/4\bZ\times \bZ/2\bZ$ and $L=\bQ(\sqrt{-5})$ is the unique subfield of $F$ with $E(L)[4]\simeq \bZ/2\bZ\times \bZ/2\bZ$. Since $L$ is imaginary and 2 is ramified in $L$, and 2 is ramified in $F'=\bQ(\sqrt{3})$, the unique real quadratic subfield of $F$. Now by theorem \ref{thm4: deg16 2}, we have $\mu_E=1$.
\end{exam}

\subsection{Examples related to further developments}\label{5.3}
In this subsection, we will present some examples that may lead to further development of the computation of $\mu$-invariants in more general cases. \par
First we will offer a counter-example of \ref{thm3: GV} when one of the elliptic curves has $\mu=1$.
\begin{exam}
Let $E: y^2+xy=x^3-520x-4225$, which is 195A4 in LMFDB, and $E': y^2+xy+y=x^3+x^2-4534140x+3499410780$, which is 2145B3. We have $E(\bQ)[4]\simeq \bZ/4\bZ\times \bZ/2\bZ\simeq E'(\bQ)[4]$ and $\bQ(E[4])=\bQ(i)=\bQ(E'[4])$. Then $E[4]\simeq E'[4]$ as both representations factor through $\Gal(\bQ(i)/\bQ)$. However, according to theorem \ref{thm4: deg2 4*2}, $\mu_E=1$ while $\mu_{E'}\geq 2$. This shows that two elliptic curves with isomorphic 4-torsion modules may not have the same $\mu$-invariant if one of them has positive $\mu$-invariant. Also, this example points out that the condition on the order of the intersection of the odd part and ramified part is necessary.
\end{exam}
According to the example above, it is natural to ask the following question.
\begin{quest}
Assume both $E$ and $E'$ have good ordinary reduction at $p=2$ and $E[8]\simeq E'[8]$ as $G_\bQ$-modules. Is it true that $\mu_E=1$ if and only if $\mu_{E'}=1$? In other words, does $E[8]$ determine $\mu$-invariant up to 1?
\end{quest}

Next we will consider an example about the fact that understanding $E[4]$ only determines the $\mu$-invariant up to 2. In other words, we are not able to tell whether $\mu$-invariant is 2, 3 or 4 just from the information of $E[4]$.
\begin{exam}
Let $E: y^2+xy+y=x^3-x^2-x-14$, this is 17A3 in LMFDB. $E': y^2+xy=x^3+605x-19750$, this is 195A8. We have $E[4]\simeq \bZ/4\bZ\simeq E'[4]$, $[\bQ(E[4]):\bQ]=4=\bQ(E'[4]):\bQ]$ and the order of the intersection of odd and ramified 4-division points are of order 4. According to theorem \ref{thm4: deg4 4}, we have $\mu_E, \mu_{E'}\geq 2$. However, according to the data in LMFDB, $E$ has analytic $\mu$-invariant 2 while the analytic $\mu$ for $E'$ is 3. If we assume Mazur's Iwasawa Main Conjecture is true, then we should expect $\mu_E=2$ and $\mu_{E'}=3$. The defect lies in the fact that according to proposition 2.10 of \cite{LY26}, $\Sel(E[4]/\bQ_\infty)^\vee$ has $\mu$-invariant $\leq 2$, so we are not able to determine the $\mu$-invariant explicitly if the 4-Selmer group has $\mu$-invariant 2.
\end{exam}
\include{Section6}
\include{Section7}
\include{Section8}
\clearpage
\addcontentsline{toc}{section}{References}
\makeatletter
\interlinepenalty=10000
\printbibliography

@incollection {Gr99,
    AUTHOR = {Greenberg, Ralph},
     TITLE = {Iwasawa theory for elliptic curves},
 BOOKTITLE = {Arithmetic theory of elliptic curves ({C}etraro, 1997)},
    SERIES = {Lecture Notes in Math.},
    VOLUME = {1716},
     PAGES = {51--144},
 PUBLISHER = {Springer, Berlin},
      YEAR = {1999}
}

@article {GV00,
    AUTHOR = {Greenberg, Ralph and Vatsal, Vinayak},
     TITLE = {On the {I}wasawa invariants of elliptic curves},
   JOURNAL = {Invent. Math.},
  FJOURNAL = {Inventiones Mathematicae},
    VOLUME = {142},
      YEAR = {2000},
    NUMBER = {1},
     PAGES = {17--63}
}

@incollection {Gr89,
    AUTHOR = {Greenberg, Ralph},
     TITLE = {Iwasawa theory for {$p$}-adic representations},
 BOOKTITLE = {Algebraic number theory},
    SERIES = {Adv. Stud. Pure Math.},
    VOLUME = {17},
     PAGES = {97--137},
 PUBLISHER = {Academic Press, Boston, MA},
      YEAR = {1989}
}

@article {Ka04,
    AUTHOR = {Kato, Kazuya},
     TITLE = {{$p$}-adic {H}odge theory and values of zeta functions of
              modular forms},
   JOURNAL = {Ast\'erisque},
  FJOURNAL = {Ast\'erisque},
    NUMBER = {295},
      YEAR = {2004},
     PAGES = {ix, 117--290}
}

@book {Wa97,
    AUTHOR = {Washington, Lawrence C.},
     TITLE = {{Introduction to cyclotomic fields}},
    SERIES = {Graduate Texts in Mathematics},
    VOLUME = {83},
   EDITION = {Second},
 PUBLISHER = {Springer-Verlag, New York},
      YEAR = {1997},
     PAGES = {xiv+487}
}

@article {Se72,
    AUTHOR = {Serre, Jean-Pierre},
     TITLE = {Propri\'et\'es galoisiennes des points d'ordre fini des
              courbes elliptiques},
   JOURNAL = {Invent. Math.},
  FJOURNAL = {Inventiones Mathematicae},
    VOLUME = {15},
      YEAR = {1972},
    NUMBER = {4},
     PAGES = {259--331}
}

@article {MSD74,
    AUTHOR = {Mazur, B. and Swinnerton-Dyer, P.},
     TITLE = {Arithmetic of {W}eil curves},
   JOURNAL = {Invent. Math.},
  FJOURNAL = {Inventiones Mathematicae},
    VOLUME = {25},
      YEAR = {1974},
     PAGES = {1--61}
}

@article {Tr05,
    AUTHOR = {Trifkovi\'c, Mak},
     TITLE = {On the vanishing of {$\mu$}-invariants of elliptic curves over
              {$\Bbb Q$}},
   JOURNAL = {Canad. J. Math.},
  FJOURNAL = {Canadian Journal of Mathematics. Journal Canadien de
              Math\'ematiques},
    VOLUME = {57},
      YEAR = {2005},
    NUMBER = {4},
     PAGES = {812--843}
}

@article {Kr99,
    AUTHOR = {Kramer, Kenneth},
     TITLE = {Elliptic curves with non-trivial {$2$}-adic {I}wasawa
              {$\mu$}-invariant},
   JOURNAL = {Acta Arith.},
  FJOURNAL = {Acta Arithmetica},
    VOLUME = {90},
      YEAR = {1999},
    NUMBER = {2},
     PAGES = {173--182}
}

@article {Ha04,
    AUTHOR = {Hachimori, Yoshitaka},
     TITLE = {On the {$\mu$}-invariants in {I}wasawa theory of elliptic
              curves},
   JOURNAL = {Japan. J. Math. (N.S.)},
  FJOURNAL = {Japanese Journal of Mathematics. New Series},
    VOLUME = {30},
      YEAR = {2004},
    NUMBER = {1},
     PAGES = {1--28}
     }

@article{LY26,
     AUTHOR = {Zichao Lin and Mulun Yin},
      TITLE = {On conjectures of {I}wasawa {$\mu$}-invariant for {S}elmer groups at {$p=2$}.},
      year = {2026},
    JOURNAL = {In preparation}
}

@article{LY26b,
     AUTHOR = {Zichao Lin and Mulun Yin},
      TITLE = {On {G}reeberg's conjecture for rational elliptic curves at {E}isenstein primes.},
      year = {2026},
    JOURNAL = {In preparation}
}
\makeatother

\end{document}